\input amstex
\input Amstex-document.sty

\pageno 511

\def\shorttitle{Weyl Manifolds and Gaussian Thermostats}
\def\shortauthor{Maciej P. Wojtkowski}

\topmatter
\title\nofrills{\boldHuge Weyl Manifolds and Gaussian Thermostats}
\endtitle

\author \Large Maciej P. Wojtkowski* \endauthor

\thanks *Department of Mathematics, University of
Arizona, Tucson, Arizona 85721, USA. E-mail: maciejw\@math.arizona.edu
\endthanks

\abstract\nofrills \centerline{\boldnormal Abstract}

\vskip 4.5mm

{\ninepoint A relation between Weyl connections and Gaussian thermostats is exposed and exploited.

\vskip 4.5mm

}
\endabstract
\endtopmatter

\document

\baselineskip 4.5mm \parindent 8mm

\specialhead \noindent \boldLARGE 1. Introduction \endspecialhead

We consider a class of mechanical dynamical systems with forcing
and a thermostating term based on the Gauss' Least Constraint
Principle for nonholonomic constraints, \cite{G99}. It was
originally proposed as a model for systems out of equilibrium,
\cite{HHP87},\, \cite{GR97},\, \cite{G99'},\, \cite{R99}.

Let us consider a  mechanical system of the form
$$
\dot{q} = v, \ \ \ \ \ \ \ \ \ \ \dot{v} = -\frac{\partial
W}{\partial q} + E,
$$
where $q,v$ from $\Bbb R^n$ are the configuration and velocity
coordinates, $W= W(q)$ is a potential function describing the
interactions in the system and $E = E(q)$ is an external field
acting on the system. The total energy $H=\frac 12 v^2 + W$ does
change because of the effect of the field $E$. We modify the
system by applying the Gauss' Principle to the constraint $H=h$ to
obtain the {\it isoenergetic} dynamical system
$$
\dot{q} = v, \ \ \ \ \ \ \ \ \ \ \dot{v} = -\frac{\partial
W}{\partial q} + E - \frac{\langle E,v \rangle}{v^2}v, \tag{1.1}
$$
with possible singularities where $v= 0$. In the special case when
$W= 0$ we have the {\it isokinetic} dynamics. The thermostating
term in the equations \thetag{1.1} is called the {\it Gaussian
thermostat}. The numerical discovery, \cite{ECM90} that the
Lyapunov spectrum, at least in the isokinetic case, has the
shifted hamiltonian symmetry raised the issue of the mathematical
nature of the equations \thetag{1.1}.

On every energy level $H=h$ the equations \thetag{1.1} define a dynamical system. In the isokinetic case the
change in $h$ is equivalent to the appropriate rescaling of time and the multiplication of the external field $E$
by a scalar.

\proclaim{Example 1.2}\endproclaim Let $\Bbb T^2$ be the flat
torus with coordinates $(x,y) \in \Bbb R^2$ and $E = (a,0)$ be the
constant vector field on $\Bbb T^2$. The Gaussian thermostat
equations on the energy level $\dot{x}^2+\dot{y}^2=1$,
$$
\ddot{x} = a\dot{y}^2 , \ \ \ddot{y} = -a\dot{x}\dot{y},
$$
can be integrated and we obtain as trajectories translations of
the curve $$ax = - \ln \cos ay$$
 or the horizontal lines.  Assuming that $E$ has an
irrational direction on $\Bbb T^2$ we obtain the following global
phase portrait for the isokinetic dynamics. In the unit tangent
bundle $S\Bbb T^2 = \Bbb T^3$ we have two invariant tori $A$ and
$R$ with minimal quasiperiodic motions, $A$ contains the unit
vectors in the direction of $E$ and it is a global attractor and
$R$ contains the unit vectors opposite to $E$ and it is a global
repellor. It can be established (\cite{W00}) that these invariant
submanifolds are normally hyperbolic so that the phase portrait is
preserved under perturbations. This example reveals that Gaussian
thermostats, even in the most restricted isokinetic case are not
in general hamiltonian with respect to any symplectic structure.
In part of the phase space they may contract phase volume and
hence have no absolutely continuous invariant measure.

The involution $I(q,v)=(q,-v)$ conjugates the forward and backward
in time dynamics, i.e., the system is reversible. Reversibility is
close to the hamiltonian property, for instance, when accompanied
by enough recurrence it can replace symplecticity in KAM
theory,\cite{Se98}.

Dettmann and Morriss, \cite {DM96}, proved the shifted symmetry of
the Lyapunov spectrum, in the case of isokinetic dynamics with a
locally potential field $E$ by exposing the locally hamiltonian
nature of the equations. For the system \thetag{1.1} with $W=0$
and $E = -\frac{\partial U}{\partial q}$ the change of variables
$$
p = e^{-U}v, \ \ \ \frac{dt}{d\tau} = e^U,
$$
brings \thetag{1.1} to the form
$$
\frac{d q}{d\tau} =\frac{\partial H}{\partial p},\ \ \ \frac{d
p}{d\tau} =-\frac{\partial H}{\partial q},\ \ \ H= \frac 12
e^{2U}p^2 = \frac 12 v^2.
$$
Under the same assumptions, Choquard, \cite{Ch97}, found a
variational principle, also involving the factor $e^U$ which in
the physically interesting examples is multivalued, thus making
the whole description only local. Liverani and Wojtkowski, \cite
{WL98}, made the observation that although the form
$$
\sum dp\wedge dq = e^{-U}\left(\sum dv\wedge dq - dU\wedge(\sum
vdq)\right),
$$
like the coordinate system $(p,q)$ is defined only locally, the
globally defined form $\omega = \sum dv\wedge dq - dU\wedge(\sum
vdq)$ can be used to develop hamiltonian-like formalism.

The three formulations above (\cite{DM96},\cite{Ch97},\cite{WL98})
apply only to isokinetic dynamics with a locally potential field
$E= -\nabla U$. In \cite{W00} a geometric setup was proposed that
covers all cases, i.e., isoenergetic and isokinetic, with a
potential vector field $E$ and nonpotential as well. We will
describe now this setup in detail.

\specialhead \noindent \boldLARGE 2. Weyl manifolds and W-flows
\endspecialhead

Let us consider a compact $n$-dimensional riemannian manifold $M$
and its tangent bundle $TM $. The metric $g$ will be also denoted
by $\langle \cdot, \cdot \rangle $. For a smooth vector field $E$
on $M$ the equations of isokinetic dynamics (the Gaussian
thermostat) on the energy level $v^2 =1$ have the coordinate free
form
$$
\frac{dq}{dt} = v,  \frac{Dv}{dt} = E- \langle E,v\rangle
v,\tag{2..1}
$$
where $\frac{D}{dt}$ denotes the covariant derivative, i.e.,
$\frac{D}{dt} = \nabla_v$ and $\nabla$ is the Levi-Civita
connection. We obtain a flow on the unit tangent bundle $SM$ of
$M$.

Let $\varphi$ be the 1-form associated with the vector field $E$,
i.e., $\varphi(\cdot) = \langle E, \cdot \rangle$. The form
$\varphi$ and the riemannian metric $g$ define a Weyl structure on
$M$, which is a linear symmetric connection $\widehat{\nabla}$
given by the formula (cf. \cite{F70})
$$
\widehat{\nabla}_XY = \nabla_XY + \varphi(Y)X + \varphi(X)Y -
\langle X , Y\rangle E,
$$
for any vector fields $X, Y$ on $M$. The Weyl structure  is
usually introduced on the basis of the conformal class of $g$
rather than $g$ itself, but in our study we fix the riemannian
metric $g$, which plays the role of a physical space. If we change
the metric $g$ to $\widetilde g = e^{-2U}g $, then the 1-form
$\varphi$ is replaced by $\widetilde\varphi = \varphi + d U$.
Hence if the vector field $E$ has a potential, i.e., $E = -\nabla
U$ then the Weyl structure coincides with the Levi-Civita
connection of the rescaled metric $\widetilde g$.

The defining property of the Weyl connection is that it is a
symmetric linear connection $\widehat{\nabla}$ such that (cf.
\cite{F70})
$$
\widehat{\nabla}_X g = - 2\varphi(X) g, \tag{2.2}
$$
for any vector field $X$ on $M$, which is equivalent to the
requirement that the parallel transport defined by the linear
connection preserves angles. We consider the geodesics of the Weyl
connection. They are given by the equations in $TM$
$$
\frac{dq}{ds} = w, \ \ \  \frac{\widehat Dw}{ds} = 0,\tag{2.3}
$$
where $\frac{\widehat D}{ds} = \widehat{\nabla}_w$. These
equations provide geodesics with a distinguished parameter $s$
defined uniquely up to scale. It follows from \thetag{2.2} and
\thetag{2.3} that $\frac{d|w|}{ds} = - \varphi(w)|w|$. Assuming
that at the initial point $q(0)$ we have $|w| = 1$ we obtain
$$
|w| = e^{-\int_{q(0)}^{q(s)}\varphi}.
$$
This formula shows that unless the form $\varphi$ is exact we
should not expect the geodesic flow  in $TM$  of a Weyl connection
to preserve any sphere bundle. We introduce the flow
$$\Phi^t: SM \to SM,$$
which we call the {\it W-flow} for the field $E$, by parametrizing
the geodesics of the Weyl connection with the arc length given by
$g$. In other words the projection of a trajectory of $\Phi^t$ to
$M$ is a geodesic of the Weyl connection, $t$ is the arc length
parameter defined by the metric $g$ and the trajectory itself is
the natural lift of the geodesic to $SM$. By direct calculation we
obtain

\proclaim{Theorem 2.4} The isoenergetic dynamics
$$
\frac{dq}{dt} = v, \ \ \   \frac{Dv}{dt} = -\nabla W +E-
\frac{\langle E,v\rangle}{v^2} v. \tag{2.5}
$$
on the energy level $\frac 12 v^2 + W =h$, reparametrized by the
arc length defines a flow on $SM$ which coincides with the W-flow
for the vector field
$$
\widetilde E = \frac{-\nabla W +E}{2(h-W)}.
$$
(We assumed implicitly that $v^2$ does not vanish on the energy
level set.) In particular in the isokinetic case on the energy
level $v^2 = 1$ we obtain that the equations \thetag{2.1} define
the W-flow for the field $E$ itself.
\endproclaim
This theorem can be interpreted as a generalization of the
Maupertuis metric,\cite{A89}. Indeed, in the case of $E=0$ we have
$$
\widetilde E = \frac{-\nabla W }{2(h-W)} = -\nabla\left(-\frac 12
\ln(h-W)\right),
$$
and hence the Weyl connection for the field $\widetilde E$ is the
Levi-Civita connection for the Maupertuis metric $(h-W)g$.

In this formulation it becomes transparent how the isoenergetic
case differs from the isokinetic. In the former case, if $E=
-\nabla U$ has a (local) potential, the vector field $\widetilde E
= \frac{-\nabla (W +U)}{2(h-W)}$ does not have a potential unless
$dW\wedge dU = 0$, i.e., unless  $W$ and $U$ are functionally
dependent. It fits well with the result of Bonetto,Cohen and Pugh
\cite{BCP99} that the isoenergetic dynamics does not in general
posses the shifted hamiltonian symmetry of the Lyapunov spectrum.
This result implies that the W-flows for nonpotential vector
fields are not in general conformally symplectic for any choice of
the conformally symplectic structure,\cite{WL98}.

\vskip 0.4cm

{\boldLARGE \noindent 3. Jacobi equations, curvature of Weyl
connec-

tions and linearizations of W-flows}

\vskip 0.4cm


We are interested in studying hyperbolicity of the dynamical
system \thetag{2.1} on $SM$. Since it was revealed to be a
modification of the geodesic flow of a connection, it is natural
to check if the Anosov theory of riemannian geodesic flows can be
extended in this direction. The first step in such a study must be
the investigation of linearized equations of \thetag{2.1}. For
riemannian geodesic flows a very useful linearization is furnished
by the Jacobi equations. The Jacobi equations are valid not only
for a Levi-Civita connection but for any symmetric connection. To
describe them let us consider a one parameter family of geodesics
of the Weyl connection, i.e., a family of solutions of
\thetag{2.3} parametrized by the parameter $u$ close to zero,
$$q(s,u),w(s,u)= \frac{dq}{ds}, |u| < \epsilon.$$
We introduce the Jacobi field
$$\xi = \frac{dq}{du} \ \ \ \text{ and } \ \ \
\widehat{\eta} = \widehat{\nabla}_\xi w.$$ The Jacobi equations
read
$$
\frac{\widehat D\xi}{ds} = \widehat{\eta}, \ \ \ \frac{\widehat
D\widehat{\eta}}{ds} = -  \widehat{R}(\xi,w)w , \tag{3.1}
$$
where for any tangent vector fields $X,Y$,
$$\widehat{R}(X,Y) = \widehat{\nabla}_X\widehat{\nabla}_Y -
\widehat{\nabla}_Y\widehat{\nabla}_X - \widehat{\nabla}_{[X,Y]}$$
 is the curvature tensor of the Weyl connection.

Let us split the vector field $\widehat{\eta} = \widehat{\eta}_0 +
\widehat{\eta}_1$ into the component $\widehat{\eta}_1$ orthogonal
to $w$ and the component $\widehat{\eta}_0$ parallel to $w$.  The
equations \thetag{3.1} can be split accordingly
$$
\frac{\widehat D\widehat{\eta}_1}{ds} = - \widehat{R}_a(\xi,w)w, \
\ \ \frac{\widehat D\widehat{\eta}_0}{ds} = -
\widehat{R}_s(\xi,w)w, \tag{3.2}
$$
where $\widehat{R}_a(X,Y)$ is the antisymmetric and
$\widehat{R}_s(X,Y)$ the symmetric part of the Weyl curvature
operator $\widehat{R}(X,Y) = \widehat{R}_a(X,Y) +
\widehat{R}_s(X,Y)$ ($\widehat{R}_s$ is called the distance
curvature and $\widehat{R}_a$ the direction curvature, cf.
\cite{F70}).

We are faced now with two tasks, to derive the linearization of
\thetag{2.1} from the equations \thetag{3.1} and to study the
curvature tensor of the Weyl connection.

Together with the family of Weyl geodesics let  us consider the
respective family of trajectories of the W-low \thetag{2.1}
$$q(t,u),v(t,u) = \frac{dq}{dt}, |u| < \epsilon.$$
We define again the Jacobi field by $\xi = \frac{dq}{du}$. Letting
$\eta = \nabla_\xi v = \nabla_v\xi$ we can consider $(\xi,\eta)$
as coordinates in the tangent space of $SM$, which is described by
$\langle v, \eta \rangle = 0$. Further we replace $\eta$ with
$\chi$, the component of $\widehat{\nabla}_\xi v$ orthogonal to
$v$. It can be calculated that
$$\chi = \eta + \langle
E,v \rangle \xi - \langle \xi,v \rangle E$$ and hence we can use
$(\xi,\chi)$ as linear coordinates in the tangent bundle of $SM$.
Note that $\langle v, \eta \rangle = 0$ is equivalent to $\langle
v, \chi \rangle = 0$ and in these new coordinates the velocity
vector field of the W-flow \thetag{2.1} is simply $(v,0)$. Now the
linearization of \thetag{2.1} can be written as
$$
\frac{\widehat D\xi}{dt} = \chi + \varphi(\xi)v, \ \ \
\frac{\widehat D\chi}{dt} = -  \widehat{R}_a(\xi,v)v
+\varphi(v)\chi. \tag{3.3}
$$
We will rewrite the equations \thetag{3.3} as ordinary linear
differential equations with time dependent coefficients. To
achieve that we need to choose frames in the tangent spaces of $M$
along the trajectory where we linearize the W-flow. The Weyl
parallel transport along a path is a conformal linear mapping and
the coefficient of dilation is equal to $e^{-\int \varphi}$. We
choose an orthonormal frame $v, e_1, \dots , e_{n-1}$ in an
initial  tangent space $T_{q_0}M$ and parallel transport it along
a trajectory of our  W-flow in the direction $v \in SM$. We obtain
the orthogonal frames which we normalize by the coefficient
$e^{\int \varphi}$ and denote them by
$v(t),e_1(t),\dots,e_{n-1}(t)$. Let $(\xi_0,\xi_1,\dots,
\xi_{n-1})$  and $(0,\chi_1,\dots, \chi_{n-1})$ be the components
of $\xi$ and $\chi$ respectively in these frames. Let further
$\widetilde \xi = (\xi_1,\dots, \xi_{n-1}) \in \Bbb R^{n-1}$ and
$\widetilde \chi = (\chi_1,\dots, \chi_{n-1}) \in \Bbb R^{n-1}$.
The equations \thetag{3.3} will read then
$$
\frac{d\xi_0}{dt} = \varphi(\xi -\langle \xi,v\rangle v), \ \ \
\frac{d\widetilde\xi}{dt} =  - \varphi(v)\widetilde \xi +
\widetilde \chi, \ \ \ \frac{ d\widetilde\chi}{dt} = -\Cal
R\widetilde\xi \tag{3.4}
$$
where the operator $\Cal R\widetilde\xi =\widehat{R}_a(\xi,v )v
\in \Bbb R^{n-1}$ (the vector $\widehat{R}_a(\xi,v)v$, being
orthogonal to $v$, is considered as an element in $\Bbb R^{n-1}$
by the expansion in the basis $e_1(t),\dots,e_{n-1}(t)$). Note
that $\widehat{R}_a(\xi,v )v$ and $\varphi(\xi -\langle
\xi,v\rangle v)$ depend $\widetilde \xi \in\Bbb R^{n-1}$ alone.

The linearized equations \thetag{3.4} for
$(\widetilde\xi,\widetilde\chi)$ differ from the Jacobi equations
in the riemannian case by the presence of the ``nonconservative''
term $- \varphi(v)\widetilde \xi$ in the first equation and the
properties of the operator $\Cal R$, which although defined
analogously in terms of the curvature tensor is not in general
symmetric. The curvature tensor can be calculated directly,
\cite{W00}, but the result is somewhat cumbersome. However the
sectional curvatures $\widehat K(\Pi)$ of the Weyl connection in
the direction of a plane $\Pi$ defined as
$$
\widehat K(\Pi) = \langle \widehat{R}_a(X,Y)Y,X\rangle,
$$
for any orthonormal basis $\{X,Y\}$ of $\Pi$, are pleasantly
transparent
$$
\widehat K(\Pi) = K(\Pi) - E_{\perp}^2 - div_\Pi E,\tag{3.5}
$$
where $K(\Pi)$ is the riemannian sectional curvature in the
direction of $\Pi$, the vector field $E_{\perp}$ is the component
of $E$ orthogonal to $\Pi$ and $div_\Pi E = \langle \nabla_{X} E,
X \rangle + \langle \nabla_{Y} E, Y \rangle$, the partial
divergence of the vector field $E$, i.e., the exponential rate of
growth of the area in the direction $\Pi$ under the flow in $M$ of
the vector field $E$.

There is a problem with this sectional curvature. It does not
depend on the Weyl connection alone, it is also effected by the
choice of the riemannian metric $g$ in the conformal class.
However the sign of sectional curvatures is well defined.

If $M$ is 2-dimensional then $E_{\perp}=0$. Moreover on a compact
manifold $M$ with a Weyl connection there is a unique metric in
the conformal class, called the Gauduchon gauge, \cite{Ga84}, such
that the vector field $E$ is divergence free. We obtain that in
dimension 2 the curvature of a Weyl connection with respect to the
Gauduchon gauge is equal to the gaussian curvature of the
Gauduchon gauge.

Let us summarize our discussion. We have a workable linearization
\thetag{3.4} of the dynamical system \thetag{2.1} and a geometric
tool (the sectional Weyl curvature \thetag{3.5}) to describe its
properties. We are ready to draw conclusions about hyperbolic
properties of the W-flows under the assumption  of negative Weyl
sectional curvature.

\specialhead \noindent \boldLARGE 4. Hyperbolic properties of
W-flows \endspecialhead

We obtain the information about hyperbolicity of W-flows studying
the quadratic form $\Cal J$ in the tangent spaces of the phase
space $SM$, defined by $\Cal J(\xi,\chi) = \langle \xi ,\,
\chi\rangle$. The form $\Cal J$ factors naturally to the quotient
bundle (the quotient by the span of the vector field \thetag{2.1},
i.e., in the $(\xi,\chi)$ coordinates the quotient by the span of
$(v,0)$). The quotient space can be represented by the subspace
$\langle \xi, v \rangle = 0$, but this subspace is not invariant
under the linearization \thetag{3.4} of the flow. The form $\Cal
J$ in the quotient space is nondegenerate and it has equal
positive and negative indices of inertia.

We take the Lie derivative of $\Cal J$ and obtain
$$
\frac{d}{dt} \Cal J   = \chi^2\ - \varphi(v)\Cal J -
\widehat{K}(\Pi)\xi^2,\tag{4.1}$$ where $\Pi$ is the plane spanned
by $v$ and $\xi$. Because of the middle term, which is absent in
the riemannian case, the negativity of the sectional curvature
does not guarantee that \thetag{4.1} is positive definite. However
it does have a weaker property that it is positive where $\Cal J =
0$. We call a flow with this property {\it strictly $\Cal
J$-separated},\cite{W01}. A strictly $\Cal J$-separated flow has a
dominated splitting, \cite{\~M84}, i.e., it has a continuous
splitting into ``weakly stable'' and ``weakly unstable'' subspaces
on which the rates of growth are uniformly separated, but which
are not necessarily decay and growth respectively. For example the
dominated splitting allows exponential growth in the ``weakly
stable'' subspace, but then all of the ``weakly unstable''
subspace grows with a bigger exponent. It needs to be stressed
that the splitting is done in the quotient space, because in
contrast to the riemannian/contact case we do not have a priori
any invariant subspaces transversal to the flow direction.

It turns out that negative sectional Weyl curvatures guarantee
even more hyperbolicity. \proclaim{Theorem 4.2} If the sectional
curvatures of the Weyl structure are negative everywhere in $M$
then the W-flow is strictly $\Cal J$-separated and hence it has
the dominated splitting into the invariant subspaces $\Cal E^+$
and $\Cal E^-$. Moreover there is uniform exponential growth of
volume on $\Cal E^+$ and uniform exponential decay of volume on
$\Cal E^-$.
\endproclaim
\proclaim{Corollary 4.3} If the sectional curvatures of the Weyl
structure are negative everywhere in $M$ then for any ergodic
invariant measure of the W-flow the largest Lyapunov exponent is
positive and the smallest Lyapunov exponent is negative.
\endproclaim
We can apply this corollary  to an individual periodic orbit and
we obtain linear instability. Moreover there are also no repelling
periodic orbits under the assumption of negative sectional Weyl
curvature.

The 2-dimensional case is special. We have \proclaim{Theorem 4.4}
For a 2-dimensional compact surface $M$ if the curvature of the
Weyl structure is negative, i.e, $\widehat K = K  - div E < 0 $ on
$M$, then the W-flow is a transitive  Anosov flow.
\endproclaim

For a locally potential vector field $div E = -\triangle U$ and we
get \proclaim{Corollary 4.5} If $K  < -\triangle U $ on a
2-dimensional surface $M$ then the W-flow is a transitive Anosov
flow.
\endproclaim

\proclaim{Corollary 4.6} If the local potential is harmonic and
the Gaussian curvature $K < 0$  on $M$ then the W-flow is a
transitive Anosov flow.
\endproclaim

We conclude that in the case of fields given by automorphic forms
on surfaces of constant negative curvature, which were studied by
Bonetto, Gentile and Mastropietro, \cite{BGM00}, the flow is
always Anosov. Further it follows from the theory of SRB measures
that if  such a flow is  Anosov then it is also automatically
dissipative, i.e., the SRB measure is singular, \cite{W00'}.

Note that in this situation we can multiply the vector field $E$
by an arbitrary scalar $\lambda$ and we still get a transitive
Anosov flow. It would be interesting to understand the asymptotics
of the SRB measure as $\lambda \to \infty$. Is the limit supported
on the union of the integral curves of $E$?  Let us stress that
this scenario differs from the perturbative conditions in
\cite{Go97}, \cite{Gr99}, \cite{W00'}, where the geodesic
curvature of the trajectories cannot be too large. Our
trajectories may have arbitrarily large geodesic curvatures and
yet they form a transitive Anosov flow.

\vskip 0.4cm

{\boldLARGE \noindent 5. Examples, extensions, open problems and
disa-

ppointments}

\vskip 0.4cm

A. In view of Theorems 4.2 and 4.4 it is natural to ask, if the negative sectional Weyl curvature is enough to
guarantee the Anosov property for the W-flow. It follows immediately from \thetag{4.1} and \thetag{3.5} that if
for every plane $\Pi$ the sectional Weyl curvatures satisfy
$$ \widehat K(\Pi) +  \frac 14 E_\Pi^2 =
K(\Pi) - div_\Pi E - E^2 + \frac 54 E_\Pi^2 < 0,
$$
then the W-flow is Anosov. We propose \proclaim{Conjecture 5.1}
There are manifolds of dimension $\geq 3$ and tangent vector
fields $E$ such that the Weyl sectional curvatures are negative
everywhere but the W-flow is not Anosov.
\endproclaim
It is also plausible that under the assumption of negative
sectional curvatures we can obtain W-flows which are nontransitive
Anosov flows, as in the examples of Franks and Williams,
\cite{FW80}.

To resolve these issues we would like to construct examples of
Weyl manifolds with negative sectional curvatures which are not
small deformations of riemannian metrics of negative sectional
curvature. In that direction we found some obstructions.

\proclaim{Proposition 5.2} There are no Weyl structures with
negative sectional curvatures in a small neighborhood of the
homogeneous Weyl structure on an n-dimensional torus (as in
Example 1.2).
\endproclaim
\proclaim{Conjecture 5.3} There are no Weyl structures of negative
sectional curvature on n-dimensional tori.
\endproclaim
It is so for $n=2$ since for 2-dimensional manifolds the Weyl
curvature is equal to the gaussian curvature of the Gauduchon
gauge.

The presence of a negative term in the formula for the Weyl
sectional curvature \thetag{3.5} gives the impression that it is
easier to find manifolds with negative Weyl curvature than with
negative riemannian curvature. The following two observations
suggest that it is not necessarily so.

If $(M_i,g_i,E_i), i = 1,2$ are two Weyl manifolds then their
cartesian product has a natural Weyl structure, and the Weyl
sectional curvature in the direction of the plane spanned by
$(E_1,0)$ and $(0,E_2)$ is either mixed (positive, negative and
zero) or always zero (iff $|E_i| =const, i=1,2$). Hence just like
in the riemannian case, product manifolds cannot have negative
sectional curvature.

Secondly, if we look for interesting homogeneous Weyl structures
we are confronted with the phenomenon that on symmetric spaces of
noncompact type, \cite{H78}, there are no homogeneous Weyl
structures at all (except for the riemannian metric itself). The
only simply connected homogeneous  riemannian  manifolds with a
compact factor and nonpositive sectional curvature are symmetric
spaces,\cite{AW76}.

\proclaim{Conjecture 5.4} The only simply connected homogeneous
Weyl manifolds with a compact factor and negative Weyl sectional
curvatures are riemannian symmetric spaces.
\endproclaim

We obtain a homogeneous Weyl manifold by taking a left invariant
metric $g$ on a Lie group and a left invariant vector field $E$.
It is known that on unimodular Lie groups the only left invariant
metrics with nonpositive
 riemannian sectional curvature are metrics with zero curvature,
\cite{AW76},\cite{M76}. In contrast, by \thetag{3.5} the
$n$-dimensional torus, $n\geq 3$, with a constant vector field $E$
has negative Weyl sectional curvature in the direction of any
plane except for the planes that contain $E$, where it vanishes.

Another instructive example is the 3-dimensional Lie group
SOL,\cite{T98}, which has mixed riemannian sectional curvatures,
but for one of the left invariant fields the Weyl sectional
curvature is nonpositive, with some negative curvature.

A weaker version of Conjecture 5.4 is \proclaim{Conjecture 5.5}
There are no left invariant Weyl structures with negative Weyl
sectional curvatures on unimodular Lie groups.
\endproclaim

B. Billiard systems are a natural extension of geodesic flows.
Similarly we can consider billard W-flows on Weyl manifolds with
boundaries by augmenting the continuous time dynamics with elastic
reflections in the boundaries. For example we can remove from a
Weyl manifold some subsets (``obstacles''). In case of convex
obstacles in a cube (or a flat torus), we obtain Sinai billiards,
which have good statistical properties.

We can introduce Weyl strict convexity of the obstacles by
requiring that the Weyl geodesics in the exterior of the obstacle
can have locally at most one point in common with the obstacle. It
can be calculated that this property has the following
infinitesimal description. Let $N$ be the field of unit vectors
orthogonal to the obstacle and pointing out of it. The riemannian
convexity of the obstacle at  a point is defined by the positive
definiteness of the riemannian shape operator $K\xi_0 =
\nabla_{\xi_0}N$, where $\xi_0$ is from the tangent subspace to
the obstacle. Similarly we introduce the operator
$$\widehat K\xi_0 = K\xi_0 +
\langle N,E\rangle \xi_0,\tag{5.6}
$$
which is the orthogonal projection of the '' Weyl shape operator''
$ \widehat{\nabla}_{\xi_0}N$ to the tangent subspace. An obstacle
is (strictly) Weyl convex if $\widehat K$ is positive (definite)
semidefinite.

Assuming that the obstacles are Weyl convex we obtain hyperbolic
properties of the billiard W-flows, \cite{W00}, in parallel with
Section 4.

Two dimensional Lorentz gas with round scatterers of radius $r$, a
constant electric field $E$ and the Gaussian thermostat is a model
of this kind. It was studied numerically by Moran and Hoover,
\cite{MH87}. Chernov et al,\cite{CELS92}, obtained exaustive
rigorous results about its SRB measures in the case of small
fields and finite horizon. We can prove the uniform hyperbolicity
of the model when $r|E| < 1$ and this inequality is sharp. Indeed,
the exponential map $F(z) = e^{|E|z}, $ takes the trajectories of
the W-flow onto the straight lines and being conformal it respects
the reflections from the boundary. Hence the mapping $F$ takes the
billiard W-flow into a billiard flow, at least locally (globally
we get billiards on the Riemann surface of the logarithm). However
the image of a disk of radius $r$ under the exponential mapping is
strictly convex if and only if $r|E| < 1$. Once the obstacles
loose convexity we readily find elliptic periodic orbits which
rules out global hyperbolicity, \cite{W00}.

C. For obscure reasons in hamiltonian systems of many particles interacting by a pair potential, which are
expected to have in general good statistical properties (and do have them in numerical experiments), hyperbolicity
in all of the phase space is rarely encountered. The notable exceptions are the Boltzmann-Sinai gas of hard
spheres,\cite{S00}, and also the one dimensional systems of falling balls, \cite{W98},\cite{W99}. In particular,
beyond the 2 dimensional examples of Knauf, \cite{K88},(which have Weyl counterparts, \cite{W00}), we do not know
of systems equivalent to geodesic flows on manifolds with  negative sectional curvatures.

For systems of particles in an external field the Gaussian
thermostat provides additional interactions and the resulting
system is not hamiltonian. Examining the simplest examples we find
that also in this case the zero and positive Weyl sectional
curvatures are common.

For noninteracting particles we get a cartesian product of Weyl
manifolds and hence we get zero sectional curvature in some
directions, as in A.

When  we consider the system of two elastic disks in the
2-dimensional torus, the cylinders that are cut out from the four
dimensional configuration space are Weyl convex, but the zero Weyl
curvature is enough to allow the presence of local simple
attractors of the type of Example 1.2, which rules out global fast
mixing and decay of correlations.

For the Lorentz gas of two noninteracting point particles in a 2
dimensional torus with round scaterrers, an external field
 and the Gaussian thermostat,
the  calculation of  \thetag{5.6} shows that the obstacles
(products of disks) in the four dimensional torus of
configurations are not Weyl convex everywhere.

These examples indicate that in the setup of Weyl geometry there
is no more freedom for the occurence of global hyperbolicity then
in the riemannian realm. The difficulty in constructing natural
examples satisfying the Chaotic Hypothesis of
Gallavotti,\cite{G01}, seems to be parallel to the scarcity of
multidimensional hamiltonian systems with strong mixing
properties, \cite{L00}.

\specialhead \noindent \boldLARGE References \endspecialhead

\widestnumber\key{CELSS93}

\ref\key{A89}\by  V.I. Arnold \book Mathematical Methods of
Classical Mechanics \publ Springer \yr 1989
\endref

\ref\key{AW76} \by R. Azencott, E.N. Wilson \paper Homogeneous
manifolds with negative curvature,I \jour Trans. AMS \yr 1976 \vol
215 \pages 323--362
\endref

\ref\key{BCP98} \by F. Bonetto, E.G.D. Cohen,  C. Pugh \paper On
the validity of the conjugate pairing rule for Lyapunov exponents
\jour J. Statist.Phys. \yr 1998 \vol 92 \pages 587--627
\endref

\ref\key{BGM00} \by F. Bonetto, G. Gentile,  V. Mastropietro
\paper Electric fields on a surface of constant negative curvature
\jour Erg. Th. Dynam. Sys. \yr 2000 \vol 20 \pages 681--696
\endref

\ref\key{CELS93} \by N.I. Chernov, G.L. Eyink, J.L. Lebowitz,
Ya.G. Sinai \paper Steady-state electric conduction in the
periodic Lorentz gas \jour Commun. Math. Phys. \vol 154 \pages
569--601\yr 1993
\endref

\ref  \key{Ch98} \by Ph. Choquard \paper Variational principles
for thermostated systems \jour Chaos \yr 1998 \vol 8 \pages
350--356
\endref

\ref\key{DM96}\by C.P. Dettmann, G.P. Morriss \paper Hamiltonian
formulation of the Gaussian isokinetic thermostat \jour Phys. Rev.
E \vol 54 \pages 2495--2500 \yr 1996
\endref

\ref  \key{FW80} \by J.M. Franks, R. Williams \paper Anomalous
Anosov flows \jour Lecture Notes in Math. \yr 1980 \vol 819 \pages
158--174
\endref

\ref  \key{F70} \by G. B. Folland \paper Weyl manifolds \jour J.
Diff. Geom. \yr 1970 \vol 4 \pages 145--153
\endref

\ref\key{G99} \by G. Gallavotti \book Statistical Mechanics \publ
Springer  \yr 1999 \endref

\ref\key{G99'} \by G. Gallavotti \paper New methods in
nonequilibrium gases and fluids \jour Open Sys. Information
Dynamics \vol 6 \yr 1999
 \pages 101--136  \endref

\ref\key{GR97} \by G. Gallavotti, D. Ruelle \paper SRB states and
nonequilibrium statistical mechanics close to equilibrium \jour
Commun. Math. Phys. \vol 190 \pages 279--285\yr 1997
\endref

\ref\key{Ga84} \by P. Gauduchon \paper La 1-forme de torsion d'une
vari\'et\'e hermitienne compacte \jour Math. Ann \vol 267 \yr 1984
\pages 495--518  \endref

\ref\key{Go97} \by N. Gouda \paper Magnetic flows of Anosov type
\jour Tohoku Math. J.\vol 49 \yr 1997 \pages 165--183  \endref

\ref\key{Gr99} \by S. Grognet \paper Flots magnetiques en courbure
negative \jour Erg. Th. Dyn. Syst. \vol 19 \yr 1999\pages 413--436
\endref

\ref\key{H78}\by  S. Helgason \book Differential Geometry, Lie
Groups and Symmetric Spaces \publ Academic Press \yr 1978
\endref

\ref  \key{K87} \by A. Knauf \paper Ergodic and topological
properties of Coulombic potentials \jour Commun. Math. Phys. \yr
1987 \vol 110 \pages 89--112
\endref

\ref\key{L00}\by  C. Liverani \paper Interacting particles \jour
Encyclopedia of Math.Sci. \vol 101 \publ Springer \yr 2000
\endref

\ref\key{\~M84}\by  R. Ma\~n\'e \paper Oseledec's theorem from the
generic viewpoint \jour Proceedings of ICM, Warsaw, 1983 \publ
PWN, Warsaw \yr 1984 \pages 1269--1276
\endref

\ref\key{Mi86} \by J. Milnor \paper Curvatures of left invariant
metrics on Lie groups \jour Adv. Math. \vol 21 \pages 293--329 \yr
1976
\endref

\ref\key{MH87} \by B. Moran, W.G. Hoover \paper Diffusion in
periodic Lorentz gas \jour J. Stat. Phys. \vol 48 \pages 709--726
\yr 1987
\endref

\ref\key{R99}\by D. Ruelle \paper Smooth dynamics and new
theoretical ideas in nonequilibrium statistical mechanics \jour J.
Stat.Phys. \vol 95 \pages 393--468 \yr 1999
\endref

\ref\key{Se98}\by M.B. Sevryuk \paper Finite-dimensional
reversible KAM theory \jour Physica D \vol 112 \pages 132--147 \yr
1998
\endref

\ref\key{S00}\by  N. Sim\'anyi \paper Hard ball systems and
semi-dispersive billiards: hypeborlicity and ergodicity \jour
Encyclopedia of Math.Sci. \pages 51--88 \vol 101 \publ Springer
\yr 2000
\endref

\ref\key{T98} \by M. Troyanov \paper L'horizon de Sol \jour
Expo.Math \vol 16 \pages 441--480 \yr 1998
\endref

\ref \key{W98} \by M.P. Wojtkowski \paper Hamiltonian systems with
linear potential and elastic constraints \jour Fundamenta
Mathematicae \vol 157 \pages 305--341  \yr 1998 \endref

\ref \key{W99} \by M.P. Wojtkowski \paper Complete hyperbolicity
in hamiltonian systems with linear potential and elastic
collisions \jour Rep. Math. Phys. \vol 44 \pages 301--312  \yr
1999
\endref

\ref \key{W00} \by M.P. Wojtkowski \paper W-flows on Weyl
manifolds and gaussian thermostats \jour J. Math. Pures. Appl \vol
79 \pages 953--974  \yr 2000\endref

\ref \key{W00'} \by M.P. Wojtkowski \paper Magnetic flows and
gaussian thermostats on manifolds of negative curvature \jour
Fundamenta Mathematicae \vol 163 \pages 177--191  \yr 2000\endref

\ref\key{W01}\by M.P. Wojtkowski \paper Monotonicity, $\Cal J$-
algebra of Potapov and Lyapunov exponents \jour Proc. Sympos. Pure
Math. \vol 169 \pages 499--521 \yr 2001
\endref

\ref\key{WL98}\by   M.P. Wojtkowski, C. Liverani \paper
Conformally symplectic dynamics and symmetry of the Lyapunov
spectrum \jour Commun. Math. Phys. \vol 194 \yr 1998 \pages 47--60
\endref

\enddocument